\documentclass[11pt]{amsart}
 \usepackage[toc,page]{appendix}
\usepackage{appendix}
\usepackage[intoc]{nomencl}
\usepackage{graphicx}
\usepackage{caption}
\usepackage{subcaption}
\usepackage{changepage}

\usepackage[backend=bibtex, style=numeric]{biblatex}
\usepackage{blkarray}

\usepackage[utf8]{inputenc}
\usepackage{listings}
\usepackage{amsmath, etoolbox}
\usepackage{tensor}
\usepackage{mathrsfs}
\usepackage{tikz-cd}
\usepackage{comment}

\usepackage{mathtools}
\numberwithin{equation}{subsection}

\usepackage{amsthm}
\usepackage{amssymb}
\usepackage{amsfonts}
\usepackage[font=it, labelfont=bf]{caption}
\usepackage{extarrows}
\usepackage{mathrsfs}
\usepackage{bm}
\usepackage[top=1.15in, bottom=1.15in, left=1.3in, right=1.3in]{geometry}

\usepackage{esint}
\usepackage[all,cmtip]{xy}
\usepackage{xcolor}
\usepackage{mathtools}
\usepackage{enumitem}

\usepackage{accents}
\usepackage{scalerel}
\usepackage{wrapfig}

\usepackage[boxsize=0.3em]{ytableau}
\usepackage[colorlinks=true]{hyperref}

\newcommand\newlink[2]{{\protect\hyperlink{#1}{\normalcolor #2}}}
\makeatletter
\newcommand\newtarget[2]{\Hy@raisedlink{\hypertarget{#1}{}}#2}
\makeatother

\usetikzlibrary{decorations.pathreplacing}

\newtheorem{thm}{{Theorem}}[section]

\newtheorem{prop}{{Proposition}}[section]
\newtheorem{defi}{{Definition}}[section]
\newtheorem{lem}{{Lemma}}[section]

\newtheorem{exa}{{Example}}[section]

\newtheorem*{Principle*}{Principle} 

\newtheorem*{thm*}{KAM Theorem}

\theoremstyle{plain}
\newtheorem{Main}{Theorem}

\newtheorem{majthm}{Theorem}

\newcommand{\SL}{\mathrm{SL}}

\newcommand{\sech}{\mathrm{sech}}
\newcommand{\R}{\mathbb{R}}

\newcommand{\Q}{\mathbb{Q}}
\newcommand{\C}{\mathbb{C}}
\newcommand{\Z}{\mathbb{Z}}
\newcommand{\N}{\mathbb{N}}

\newcommand{\T}{\mathbb{T}}

\newcommand{\cst}{\mathrm{cst}}

\newcommand{\Hess}{\mathrm{Hess}}

\makeatletter
\newcommand{\tpitchfork}{%
  \vbox{
    \baselineskip\z@skip
    \lineskip-.52ex
    \lineskiplimit\maxdimen
    \m@th
    \ialign{##\crcr\hidewidth\smash{$-$}\hidewidth\crcr$\pitchfork$\crcr}
  }%
}
\makeatother

\def\cH{\mathcal H}

\def\cP{\mathcal P}

\def\cU{\mathcal U}
\def\cV{\mathcal V}

\def\e{\varepsilon}

\def\ph{\varphi}
\def\ti{\tilde}

\def\dt{\Delta}
\def\vs{\varsigma}

\def\be{\begin{equation}}
\def\ee{\end{equation}}

\def\bmat{\begin{pmatrix}}
\def\emat{\end{pmatrix}}

\def\bsm{\begin{smallmatrix}}
\def\esm{\end{smallmatrix}}

\newcommand{\black}[1]{\color{black}}

\newcommand{\Yi}[1]{{\color{magenta}{#1}}}

\newcommand{\two}{\mbox{\scaleto{(2)}{7pt}}}

\newcommand\B{\newlink{def:Box}{\mathcal{D}}}

\date{}

\author{Vadim Kaloshin}
\address{Institute of Science and Technology Austria, Klosterneuburg, Austria}
\email{vadim.kaloshin@gmail.com}

\author{Illya Koval}
\address{Institute of Science and Technology Austria, Klosterneuburg, Austria}
\email{illyakoval2001@gmail.com}

\author{Yi Pan}
\address{Centro di Ricerca Matematica Ennio De Giorgi, Scuola Normale Superiore, Pisa, Italy}
\email{yi.pan@sns.it}

\title[Genericity of meandering tori]{Isoenergetic degeneracy generically creates meandering invariant tori}

\begin{document}
\graphicspath{ {./figures/} }
\maketitle

\begin{abstract}
Consider the set of isoenergetically degenerate integrable Hamiltonians with two degrees of freedom. We show that a cusp-generic perturbation of a generic Hamiltonian in this set gives rise to meandering invariant tori -- embedded Lagrangian tori which are not graphs. Moreover, an exponentially dense subset of perturbations admits higher order meandering tori, of all orders from two to infinity. These infinite order meanders have an endless nested structure.
\end{abstract}

\section{Introduction}
\begin{wrapfigure}{r}{0.5\textwidth}
  \centering
  \includegraphics[width=0.48\textwidth]{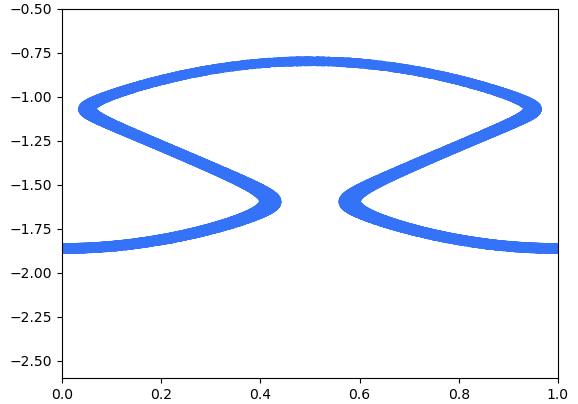}
  \caption{Poincar\'e section of $H=h+p$ at $H = -0.45$. }
  	\label{fig.meander1d1}
\end{wrapfigure}
By KAM theory, when an integrable Hamiltonian $h(I)$ satisfies a suitable non-degeneracy condition, small perturbations  preserve many invariant Lagrangian tori, which remain nearly flat. More precisely, they are Lagrangian graphs of the form 
$\{(\ph,I_0+\nabla u(\ph)):\ph\in\T^n\}$, 
where $u:\T^n\to \R$ is nearly constant. When the non-degeneracy condition fails, however, more interesting invariant tori may appear. 
For example, $h(I) = I_1^2 + I_2^3$ is isoenergetically degenerate. Taking $p(\ph)=0.04 \cos\left[2\pi(2\varphi_1 - \varphi_2) \right]$, the Hamiltonian $h+p$ admits a meander -- an embedded invariant Lagrangian torus which is not a graph, visible in the Poincar\'e section in Figure~\ref{fig.meander1d1}.


As discussed in Section~\ref{sect.ex}, meanders have been observed in many specific Hamiltonian systems, mostly by numerical methods. We provide the first rigorous result on the existence of meanders under generic deformations: for a generic isoenergetically degenerate $h(I)$ and a cusp-generic perturbation $p(\ph,I)$,  the Hamiltonian $h+p$ admits an abundance of non-graphical invariant Lagrangian tori of meandering shape. We now explain this result in detail.

\smallskip 
Let $D \subset \mathbb R^2$ be a bounded open set on the plane, $r>0$, $s>0$.
We define $V_rD=\{I\in\C^2:\|I-D\|<r\}$ to be the complex $r$-neighborhood of $D$, $W_s\T^2=\{\ph\in\T^2+i\R^2:|\Im\ph|<s\} $ to be a complex strip around the torus. We denote the set of analytic functions $h(I)$ defined on $V_rD$ by $C^\omega_r(D)$.
A function $h(I) \in C^\omega_r(D)$ is called {\it isoenergetically degenerate}\footnote{We refer to Proposition~\ref{gene.cond.main} for more discussion on the condition.} at 
$I_* \in D$ if 
\begin{equation}
\label{def.i.d} 
\nabla h(I_*)^\perp\,\Hess\: h(I_*)\,\left(\nabla h(I_*)^\perp\right)^\top=0.
\end{equation}
Denote the canonical projections by 
$\pi_I:\mathbb T^2 \times D \to D$ and $\pi_\varphi:\mathbb T^2 \times D \to \mathbb T^2$. We call the image of an analytic invertible map 
$\Phi:\mathbb T^2 \to \mathbb T^2 \times D$  an analytic embedded torus.

\begin{wrapfigure}{r}{0.4\textwidth}
  \centering
  \includegraphics[width=0.37\textwidth]{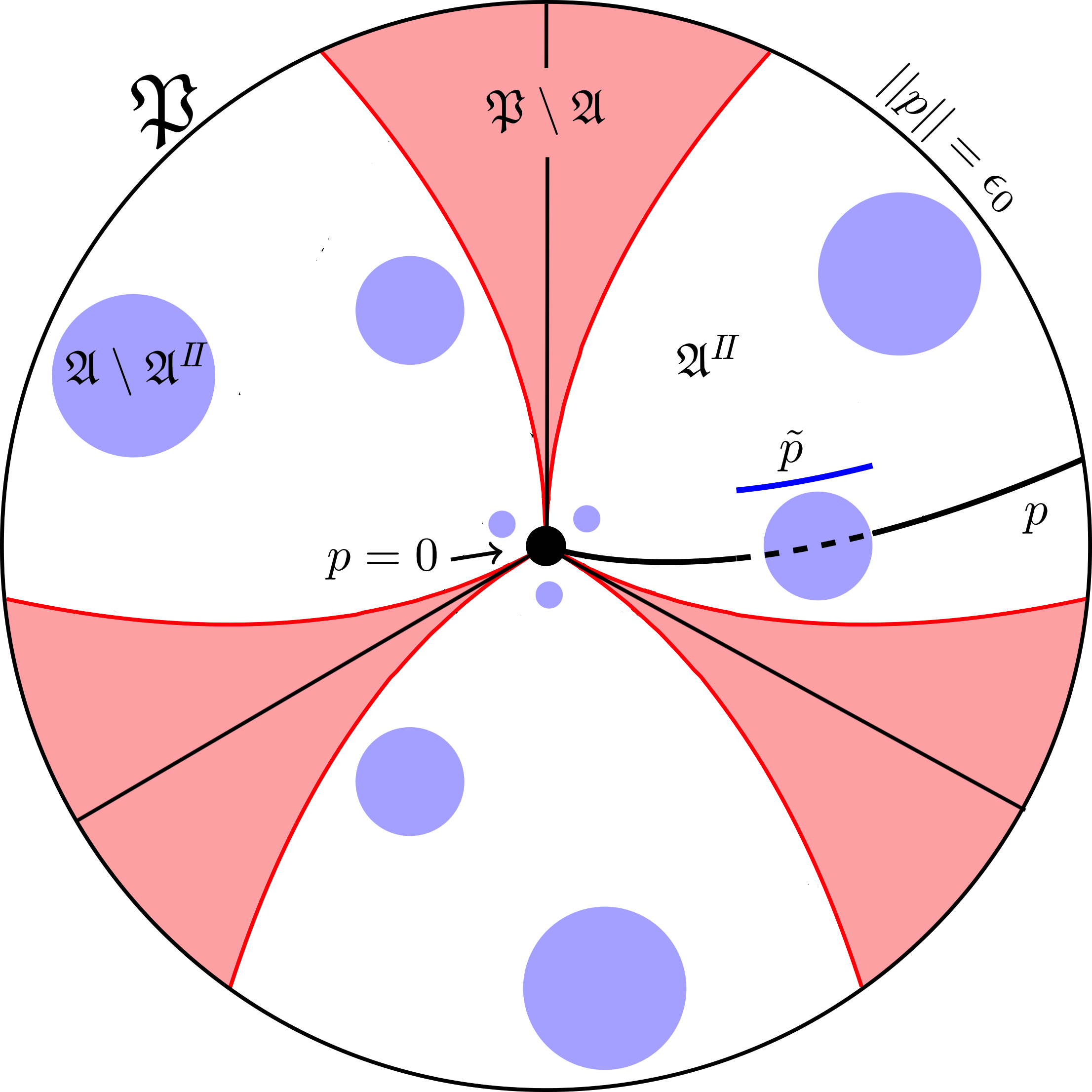}
  \caption{A cusp-generic set}
  	\label{fig.deformationspace}
\end{wrapfigure}

\begin{Main}
 \label{thm.A}
Let $h(I) \in C^\omega_r(D)$ be a generic isoenergetically degenerate Hamiltonian, $p(\varphi,I)\in C^\omega(W_s\mathbb{T}^2\times V_rD)$ be a cusp-generic\footnote{Cusp-genericity is illustrated in Figure~\ref{fig.deformationspace} and discussed in detail in Section~\ref{sect.state}.} perturbation of norm $\varepsilon$\footnote{We use the supremum norm $\|p\|_{r,s}$ on the region $W_s\T^2\times V_rD$.}.  Then $h+p$ admits an analytic embedded invariant Lagrangian torus $\mathcal{MT}$ such that:
 \begin{itemize}
  \item the $I$-projection $\pi_I(\mathcal{MT})$ lies in an $O(\varepsilon^{1/3})$-neighborhood of some $I_*$\footnote{Namely, the point at which $h$ is isoenergetically degenerate.},
 \item the $\varphi$-projection $\pi_\varphi(\mathcal{MT})$, outside two strips of width $O(\varepsilon^{3/4})$ in $\mathbb{T}^2$, is a three-to-one cover.
\end{itemize}
\end{Main}

\begin{figure}
	\includegraphics[width=7cm]{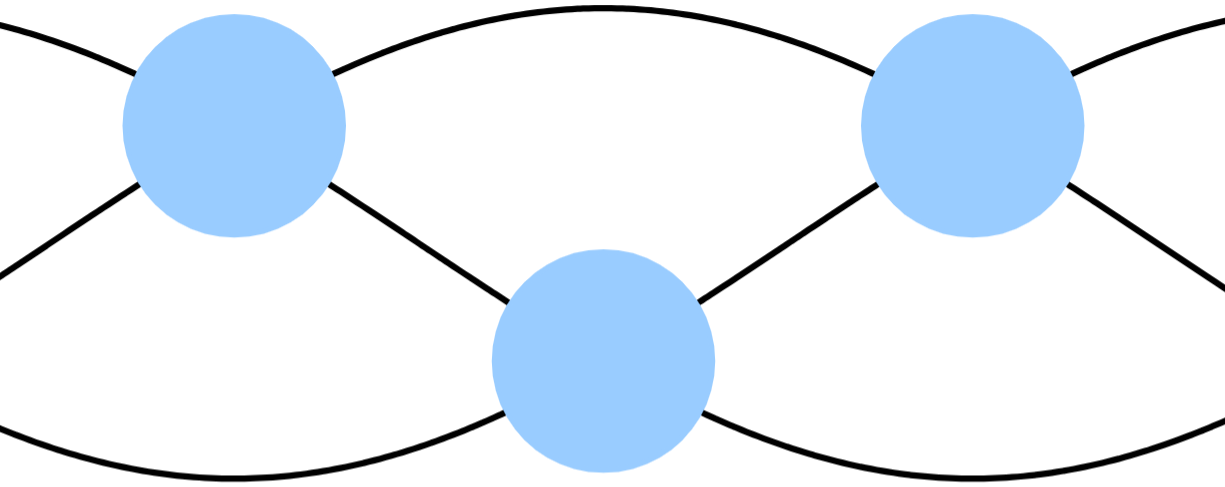}
	\centering
	\caption{A meander in $\{\ph_2=0\}$ section, on $(\varphi_2, I_2)$ plane.
}
	\label{fig.meanderdef}
\end{figure}

We call $\mathcal{MT}_1=\mathcal{MT}$ a \textbf{first order meander}. The shape of this meandering torus is shown in Figure~\ref{fig.meanderdef}. Let $(\theta_{(1)}, J_{(1)})\in\mathbb T^2 \times D_{(1)}$ be the local normal coordinates 
near ${\mathcal{MT}}.$ In particular, 
${\mathcal{MT}}=\{J_{(1)}=0\}$. Denote the canonical projections by 
$\pi_J$ and $\pi_\theta$. We have the following result.

 \begin{Main}
	\label{thm.B}
Under the assumptions of Theorem~\ref{thm.A}, for any $\eta>0$ small enough, there exists an analytic Hamiltonian $\tilde H=h + \tilde p$ with 
\[\|p -\tilde p \|_{r, s} < \exp(-1/\varepsilon),\]
such that:
\begin{itemize}
  \item $\tilde H$ has a first order meander ${\mathcal{MT}}_1$,
  \item $\tilde H$ admits an analytic 
embedded invariant Lagrangian torus ${\mathcal {MT}}_2$ such that
\begin{itemize}
  \item the $J_{(1)}$-projection $\pi_{J_{(1)}}({\mathcal{MT}}_2)$ lies in an $O(\e_{\two})$-neighborhood of some $J^*_1$,
  \item the $\theta_{(1)}$-projection $\pi_{\theta_{(1)}}({\mathcal{MT}}_2)$, outside several strips 
of width $O(\e_{\two})$ in $\mathbb T^2$, is a three-to-one cover, where $\e_{\two} = \exp(-\e^{-\eta/6})$.
\end{itemize}
\end{itemize}
\end{Main}

\begin{figure}[ht!]
	\includegraphics[width=9cm]{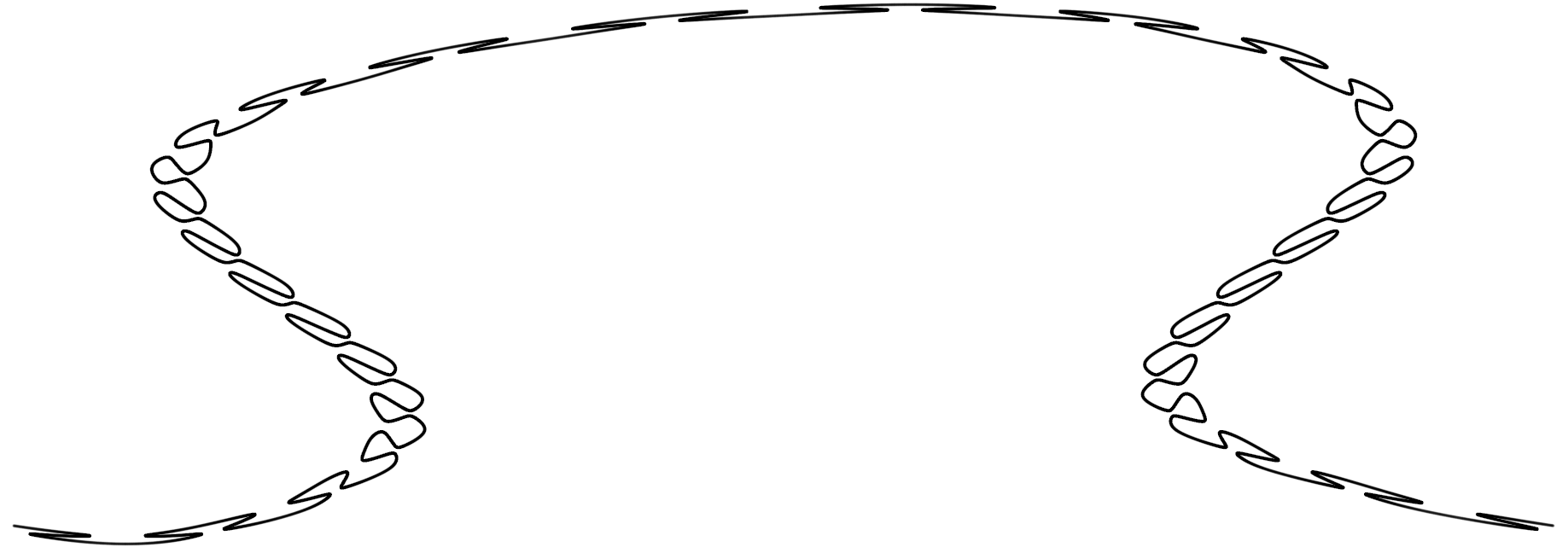}
	\centering
	\caption{Second order meandering torus $\mathcal{MT}_2$ on a section.
	}
	\label{fig.meandersec}
\end{figure}

We call $\mathcal{MT}_2$ {a \textbf{second order} meander}, an example of which is shown in Figure \ref{fig.meandersec}.  Suppose there is a countable collection of analytic embedded tori $\{\mathcal{MT}_n\}_{n\in\N^*}$. For each $n\in\N^*$, denote by $(\theta_{(n)}, J_{(n)})\in\mathbb T^2 \times D_{(n)}$ the local normal coordinates near ${\mathcal{MT}_n}$. In particular, ${\mathcal{MT}_n}=\{J_{(n)}=0\}$ in 
these coordinates. For $n \ge 2$, the auxiliary parameters $\varepsilon_{(n+1)} = \exp\left(-\varepsilon_{(n)}^{-1/6}\right)$ are defined recursively.

\begin{Main}
\label{thm.C}
Under the assumptions of Theorem~\ref{thm.A},  there exists an analytic Hamiltonian $\overline H=h + \overline p$ with  $\| p - \overline p \|_{r,s} < \exp(-1/\varepsilon)$, such that for any $n\in\N^*$:
\begin{itemize}
  \item $\overline H$ has an analytic embedded invariant Lagrangian torus $\mathcal{MT}_k$ with $(J_{(k)},\theta_{(k)})$ 
local normal coordinates, for any $k\in[1,n)\cap\N^*$, 
  \item $\overline H$ has an order $n$ analytic embedded 
invariant Lagrangian torus $\mathcal{MT}_n$ such that for any $k\in(1,n]\cap\N^*$
\begin{itemize}
  \item the $J_{(k-1)}$-projection $\pi_{J_{(k-1)}}({\mathcal{MT}}_k)$
belongs to $O(\e_{(k)})$  of some $J_{(k-1)}^*$,
  \item the $\theta_{(k-1)}$-projection $\pi_{\theta_{(k-1)}}({\mathcal{MT}}_k)$, outside of 
several strips of width $O(\e_{(k)})$ in  $\mathbb T^2$, 
is a three-to-one cover.
\end{itemize}
\end{itemize} 

The Hausdorff limit $\mathcal{MT}_{\infty}$ of $\mathcal{MT}_n$ exists. 
Moreover, $\mathcal{MT}_{\infty}$ in Poincar\'e section is an $\overline H$ invariant topological curve.
\end{Main}

We call $\mathcal{MT}_{\infty}$ an \textbf{infinite order meander} or a \textbf{nopal}. These special infinitely branching meanders were long thought to exist, with Rafael de la Llave comparing them to cacti.

If $\partial h_{I_j} (I_*)\ne  0$ for $j=1$ or $2$, then Theorems~\ref{thm.A}--\ref{thm.C}
hold for the Poincaré map on the section $\{\varphi_j=0 \}$, with meandering tori replaced by meandering curves.

Figure \ref{fig.deformationspace} is an illustration of cusp-genericity. The complement of red region in the ball $\mathfrak{P}$ is cusp-generic, which provides first order meander. The blue balls are of radius $O(\exp(-1/\varepsilon))$. In the complement of the red and the blue, a second-order or a higher-order meander exists. We refer to Definition~\ref{def.rnd} and Proposition~\ref{prop.rnd} for the precise conditions on $p$ and their cusp-genericity.

\subsection{Creation of meanders through separatrix reconnection }
\label{sect.brief.pf}
While the results in this paper are formally stated for Hamiltonians, we anticipate that similar phenomena occur for area-preserving maps. To appeal to a broader interest, we explain the creation of meanders for maps, demonstrating that the core ideas remain consistent across both settings.

Consider a standard area-preserving nontwist\footnote{For an area-preserving map $f$ on $\T\times\R$, let $F=(F_1,F_2)$ be a lift to the universal cover $\R\times\R$. We say $f$ is a {\it twist map} if $\partial_I F_1>0$. For example, replacing the definition of $\bar{\theta}$ in \eqref{eq.stnontw} by $ \bar{\theta} = \theta + 2\pi\omega_\varepsilon +\bar{I}$ yields a standard area-preserving twist map.} map.
\begin{equation}
                \begin{cases}
                    \bar{\theta} = \theta + 2\pi\omega_\varepsilon +\bar{I}^2\ (\text{mod }2\pi)\\
                    \bar{I} = I + \varepsilon U(\theta)
                \end{cases}
            \label{eq.stnontw}
            \end{equation}
            \medskip           
For example, let $\omega_\varepsilon=1/2 - \varepsilon^{2/3}$. If $U\equiv 0$, then for any small positive $\varepsilon$ there are two 
families of periodic points of period $2$ given by $\{I= \pm \varepsilon^{1/3}\}$. In other words, 
there is a break in the twist \footnote{A different way to characterize the break of twist is 
that for certain rotation numbers, associated invariant curves are not unique.}. 

Following the example in \cite{Si98}, let
$$
U(\theta)=\frac{\sin  \theta}{1+\beta \cos \theta}=
\sin \theta \sum_{k\ge 0} (-\beta)^k \cos^k \theta,\quad \beta=1/2.
$$ 
For a positive integer $q>1$ 
let $[U]_q(\theta)=\sum_{k=0}^{q-1} U(\theta +k/q)$ be the $q$-average of $U$. The map satisfies the following:
\begin{itemize} 
\item The averaging theory shows that after 
a coordinate change $O(\varepsilon)$-close to the identity, the map takes the form 
\begin{equation*}
                \begin{cases}
                    \bar{\theta} = \theta + 2\pi\omega_\varepsilon +\bar{I}^2+O(0.01)\ (\text{mod }2\pi)\\
                    \bar{I} = I + \varepsilon(0.28\sin 2\theta + 0.02\sin 4\theta +O(0.01)).
                 \end{cases}
\end{equation*}
Then for small $\varepsilon>0$, there are two saddle periodic orbits of period 
$2$, denoted by $p^\pm$, where $p^+$ is close to $\{\pi k,\ k=0,1\}\times \{I= \varepsilon^{1/3}\}$, $p^-$ is close to $\{\pi (2k+1)/2,\ k=0,1\}\times \{I= -\varepsilon^{1/3}\}$.

\begin{figure}[ht!]
\begin{subfigure}{.3\textwidth}
  \centering
			         \includegraphics[width=4cm]{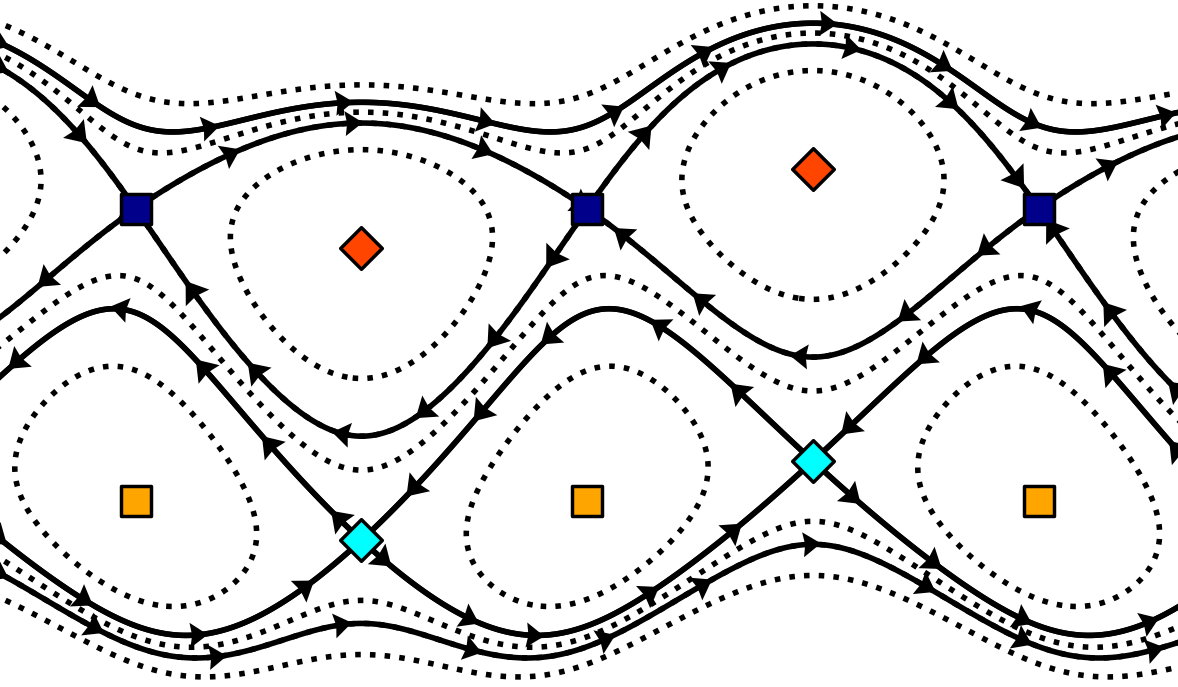}
\end{subfigure}%
\begin{subfigure}{.3\textwidth}
  \centering
  \includegraphics[width=4cm]{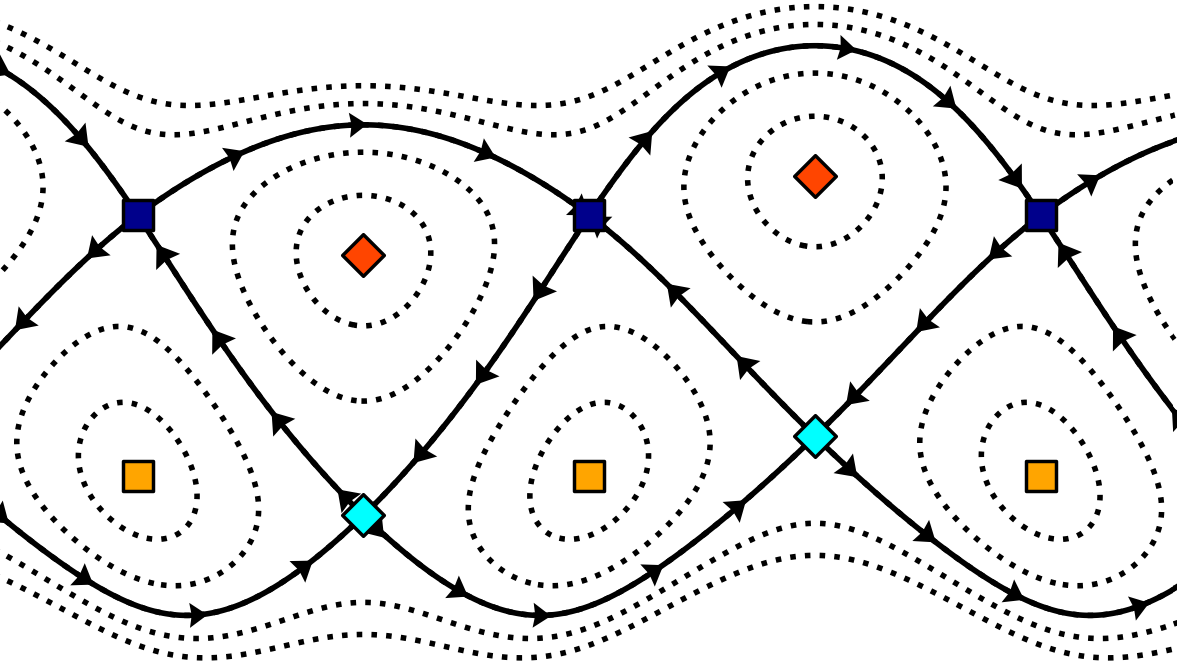}
\end{subfigure}
\begin{subfigure}{.3\textwidth}
  \centering
  \includegraphics[width=4cm]{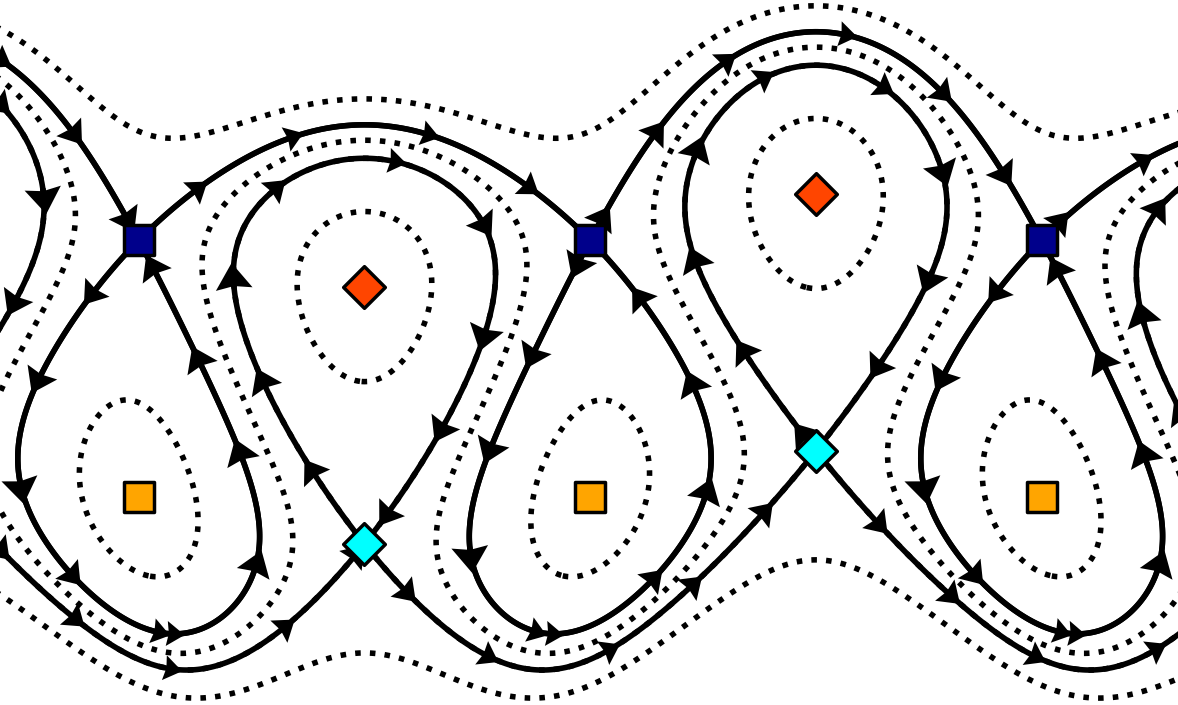}
\end{subfigure}
\caption{The dynamics of an integrable approximation of the map \eqref{eq.stnontw} as $\omega$ changes.}
\label{fig.reconnection}
\end{figure}

\item There is an almost conserved first integral 
$$H_{\omega}(\theta,I)=(\omega_\varepsilon - 1/2) I + \frac{I^3}{3} - \varepsilon\int_0^{\theta} [U]_2(\vartheta) d\vartheta.$$
For some $c^*$ and $\omega_c = 1/2 - c^* \varepsilon^{2/3}$, the level sets 
$h^+=H_\omega(p^+)$ and $h^-=H_\omega(p^-)$ coincide, which leads to a separatrix reconnection for the averaged system, see the middle of Figure \ref{fig.reconnection}.

\item  One can show that by changing $\omega \ne \omega_c$, a saddle separatix breaks and 
the averaged system has a family of meanders given by the level set $\{H=h\}$  for 
$h\in (h^-,h^+)$,  see the left and right of Figure \ref{fig.reconnection}. 
\medskip 

\item Notice that period of motion on the level set $\{H=h\}$, denoted by $T(h)$, tends to 
infinity both as $h\to h^\pm$. One can show that for small $\varepsilon>0$ there is a single 
minimum $T(h^*)=\min_{h\in (h^-,h^+)} T(h)$.
\medskip 
\item Near the level set $\{H=h\}$  for 
$h\in (h^-,h^+)$, the original system has invariant meander curves.
\medskip 

\item Break of twist is equivalent to the break of monotonicity of $T(h)$. This allows us to iterate the construction. 
\end{itemize}



\section{Precise statement of results}
\label{sect.state}
To state Theorems \ref{thm.A}-\ref{thm.C} precisely, we require the following conditions. Define
\begin{equation}
	\gamma(I):=\nabla h(I)^\perp\,\Hess\: h(I)\,\left(\nabla h(I)^\perp\right)^\top.
\end{equation}
\begin{defi}
\label{def.gid}
	An integrable Hamiltonian $h(I)$ of two degree of freedom is called \textup{generic isoenergetically degenerate} at $I_*$ if the following conditions hold:
	\begin{align}
		&\gamma(I_*)=0, \label{ConditionA}\tag{a}\\
		&\nabla h(I_*)\neq(0,0), \label{ConditionB}\tag{b}\\
		&\nabla \gamma(I_*)\cdot\nabla h(I_*)^\perp\neq0, \label{ConditionC}\tag{c}\\
		&\det \Hess\: h(I_*) < 0, \label{ConditionD}\tag{d}\\
		&\nabla h(I_*)\cdot j=0,\quad \text{for some }j\in\Z^2\setminus \{0\}.\label{ConditionE}\tag{e}
	\end{align}
\end{defi}

Among these, Conditions \eqref{ConditionA}--\eqref{ConditionB} are the most restrictive, while the remaining ones hold fairly generally. We denote by $\mathcal{ID}(D)$ the set of analytic Hamiltonians on $D$ that satisfy Conditions \eqref{ConditionA}--\eqref{ConditionB} for some $I_* \in D$.  Condition \eqref{ConditionA} is precisely the isoenergetic degeneracy we defined in \eqref{def.i.d}, which can be rewritten as
\[	\det \begin{pmatrix}
		h_{11} & h_{12} & h_1 \\
		h_{12} & h_{22} & h_2 \\
		h_1 & h_2 & 0
	\end{pmatrix}(I_*) = 0.
\]
This condition fails on an open set of Hamiltonians. For instance, it is never satisfied if $h$ is Tonelli (convex in $I$). Nevertheless, the following results show that the conditions are broadly satisfied.

\begin{prop}
	\label{gene.cond.main}
$\mathcal{ID}(D)$ is nonempty and contains an open set. Moreover, there exists an open dense subset $\cH \subset \mathcal{ID}(D)$ of generic isoenergetically degenerate Hamiltonians.
\end{prop}

\begin{lem}
	For any analytic $h$ on $\R^2$, there exists an open set of $(a,b)\in\R^2$, such that $\ti h(I)+aI_1^2+bI_2^2\in\mathcal{ID}(\mathbb{R}^2)$.
	\label{lem.big.h0}
\end{lem}

We also have conditions on deformations.
\begin{defi}
\label{def.rnd}
	Let $\vs>0$. A function $g(\xi)\in C^2(\T)$ on $\T$ is called \textup{$\vs$-non-degenerate} if there exists a unique maximum and minimum of $g(\cdot)$, denoted by $\xi_+$ and $\xi_-$ respectively, satisfying
\begin{align}
		&g(\xi_+) - g(\xi_-) > \vs, \label{ConditionF}\tag{f} \\
		&\left|g''(\xi_\pm)\right| > \vs ,\label{ConditionG} \tag{g}\\
		&\left|g(\xi) - g(\xi_\pm) \right| \ge \vs\|\xi - \xi_\pm\|^2, \quad \xi\in\T.
		\label{ConditionH}\tag{h}
	\end{align}
Let $j=(j_1,j_2)\in\Z^2\backslash\{0\}$, $\eta>0$. A deformation $p(\varphi, I)\in C^2(\T^2\times\R^2)$ is called \textup{$j$-resonant $\eta$-non-degenerate at $I_*$} if its resonant part $g(\xi)$, defined by
	$$g(\xi)=\int_0^1p(r_2\xi-j_2\ph_2,-r_1\xi+j_1\ph_2, I_*)\,d\ph_2,$$
is $\|p\|_{r, s}^{1+\eta}$-non-degenerate, where $r_1,r_2\in\Z$ satisfy that $r_2j_1-r_1j_2=1$ (existence given by B\'ezout's identity).
\end{defi}

We point out that Conditions~\eqref{ConditionF}-\eqref{ConditionH} are weaker if $\vs$ is smaller.

\begin{prop}
\label{prop.rnd}
Let $j\in\Z^2\backslash\{0\}$, $0 < \eta < 1$.
The set 
 of $j$-resonant $\eta$-non-degenerate deformations $p\in C^\omega(W_s\T^2\times V_rD)$ is \emph{cusp-generic} in the following sense:
 \begin{itemize}
   \item there is an open dense set $\cU\subset\{p(\ph,I) \in C^\omega( W_s\T^2\times V_rD), \|p\|_{r,s}\leq1\}$,
   \item there is a nonnegative continuous function $\e_0:C^\omega( W_s\T^2\times V_rD)\to\R_{\geq0}$ with $\e_0|_{\cU}>0$,
   \item there is a cusp set $\cV=\cup_{p\in \cU} \cup_{0<\e<\e_0(p)}\{\e p\}$,
 \end{itemize}
   so that any function in $\cV$ is a $j$-resonant $\eta$-non-degenerate deformation.
\end{prop}

By choosing a proper function $\e_0$, elements in the cusp set $\cV$ can have a uniform bound and be small enough. That is how we require the smallness of $\|p\|_{r,s}$ implicitly in Theorems~\ref{thm.A}-\ref{thm.C}.

\smallskip
With the above conditions in place, we now give a precise formulation of Theorems~\ref{thm.A}--\ref{thm.C}.

\begin{Main}
	\label{thm.D}
	Let $h(I)\in C^\omega_r(D)$, $p(\varphi,I)\in C^\omega(W_s\T^2\times V_rD)$, $\eta>0$ small enough. Assume that $h(I)$ is 
	generic isoenergetically degenerate at some $I_*\in D$ with 
	$j\in \mathbb{Z}^2\backslash\{0\}$, and $p(\varphi, I)$ is $j$-resonant $\eta$-non-degenerate at $I_*$ with norm $\e$ small enough. Then 
\begin{enumerate}
\item $H =h+p$ has a meandering invariant Lagrangian torus close to $I = I_*$,
\item there exists $\tilde{p}(\varphi, I)\in C^\omega(W_s\T^2\times V_rD)$  satisfying
	\begin{equation}
		\|\tilde{p}- p\|_{r,s} \le \exp\left(-1/\e\right),
	\label{eq.thmB}
	\end{equation}
such that $\ti H = h+ \tilde{p}$ has a second order meander, 
\item there exists $\overline{p}(\varphi, I)\in C^\omega(W_s\T^2\times V_rD)$ satisfying \eqref{eq.thmB}, such that $\overline H = h+ \overline{p}$ has meanders of arbitrary order $n \ge 1$, as well as a meander of infinite order.
\end{enumerate}
\end{Main}

\section{Some history and examples}
\label{sect.ex}

\subsection{A brief historical survey}
The study of meanders has a rich history. It dates back to 1956: in \cite{S56} meandering invariant tori have been found for Hamiltonian with one degree of freedom. Through the decades, many other numerical examples have been studied. We show some of them in the next subsection. 

For an area-preserving twist map on $\T\times\R$, Birkhoff \cite{B22} proved that its essential invariant curves are graphs of Lipschitz functions. Furthermore, invariant curves always exist in a nearly-integrable case,
 as shown by Moser \cite{M62}.
It is a natural question to ask what happens when the twist condition is violated.

In this nontwist setting, Herman \cite{H92} proved the existence of meandering invariant tori for some symplectic diffeomorphisms. 
Later, Sim\'o \cite{Si98} studied invaraint curves of general nontwist maps. In fact, he proposed an idea to show the existence of meanders for a generic perturbation of a nontwist area-preserving map. Furthermore, he provided a numerical example of a higher order meander. These insights and examples serve as the starting point for our study. Although our current results are established in the Hamiltonian setting, we expect similar results for area-preserving nontwist maps.

As for Hamiltonian systems, the Birkhoff theorem was extended by Arnaud \cite{Ar10} to Tonelli Hamiltonians on compact and connected manifolds, showing that zero section isotopic invariant Lagrangian submanifolds are always Lagrangian graphs. Here, Lagrangian assumption is important -- as without it Arnaud \cite{Ar14} constructed a counterexample by essentially gluing smaller dimensional invariant tori together.  Moving beyond Tonelli Hamiltonians, we work with isoenergetically degenerate Hamiltonian systems, constructing meandering invariant tori that are Lagrangian and isotopic to graphs.

More recently, González-Enríquez, Haro and De la Llave \cite{GHL13, GHL22} studied
nontwist systems and creation of meanders using the methods of bifurcation theory. 
Besides, non-essential invariant curves of meandering shape
appear to be common even for twist maps. Whenever a generic map has an elliptic periodic point whose eigenvalue passes through $\exp(2\pi i/3)$, a non-twist invariant torus arises around the periodic orbit \cite{DMS00}. In particular, this mechanism gives rise to invariant curves of meandering shape for the area-preserving Hénon map \cite{DMS00}.

\subsection{Possible applications}
Theorem~\ref{thm.D} applies to general Hamiltonian systems. It is natural to ask for concrete examples. Numerous examples of isoenergetically degenerate Hamiltonians are given in \cite{nontwist}, a selection of which we study below. Many of these examples are not generically isoenergetically degenerate, so Theorem~\ref{thm.D} does not apply directly. Moreover, the perturbation construction for higher-order meanders requires adaptation to each specific example. Nevertheless, we conjecture that suitably modified versions of Theorem~\ref{thm.D} apply to examples below.
\begin{exa}[\cite{Rossby}, Rossby wave]
	$H = -\tanh I + \sech^2 I \sum_{i = 1}^2 \e_i \cos k_i(\varphi - c_i t)$.
\end{exa}

This example models the atmosphere of the planet and it may explain the meanders we frequently see on weather maps. Here, $I$ represents the latitude, with $I = 0$ being some fixed latitude line. The longitude is represented by $\varphi$. At $I = 0$ the Hamiltonian is non-isoenergetic.  

\begin{exa}[\cite{K68}, Kepler problem]
	$H = -\frac{\mu^2}{2I_L^2} + p(\varphi, I)$.
\end{exa}
The Kepler problem in celestial mechanics is one of the most classical integrable systems. In Delaunay variables, the Hamiltonian depends only on one action $I_L$, making it always isoenergetically degenerate. Integrability can be destroyed in various ways. For instance, in \cite{K68} the planet is assumed to be slightly oblate rather than perfectly spherical, modeling the actual shape of the Earth. In the related planar Hill problem, Sim\'o-Stuchi~\cite{SS00} numerically observed meandering invariant torus.

\subsection{Some examples of Hamiltonians}
We now examine several examples of Hamiltonians and verify whether they satisfy Conditions~\eqref{ConditionA}--\eqref{ConditionE} for generic isoenergetic degeneracy.
With respect to the conditions, we define 
\[\Gamma=\{I~|~\gamma(I)=0\},\quad \Pi_j=\{I~|~\nabla h(I)\cdot j=0\},~~j\in\Z^2.\]

\begin{exa}
	$h(I) = I_1^2 + I_2^2$.
\end{exa}

This is a standard Tonelli Hamiltonian. The Hessian $\Hess \:  h(I)$ is positive definite everywhere, so $\gamma(I_*)$ is always positive unless $\nabla h(I_*) = 0$, which violates Condition~\eqref{ConditionB}. Hence $h \notin \mathcal{ID}(\R^2)$. Furthermore, this convexity persists under small perturbations.

\begin{exa}
	$h(I) = I_1^4 + I_2^4 + I_1^2 + I_2^2 - 4I_1 I_2 + 10I_1$.
	\label{example.nonablnabl}
\end{exa}
This example is not Tonelli since the Hessian is indefinite in a neighborhood $I = 0$. In this region, the quadratic form $\Hess\: h(I)$ has two null directions. However, none of these directions are ever aligned with the gradient, which is dominated by the term $10I_1$. This forces $\Gamma$ to remain empty despite the indefiniteness of $\Hess\: h$, and this obstruction persists under small deformations. This is an example of a quasi-convex Hamiltonian.

\begin{figure}[ht!]
	\includegraphics[width=10cm]{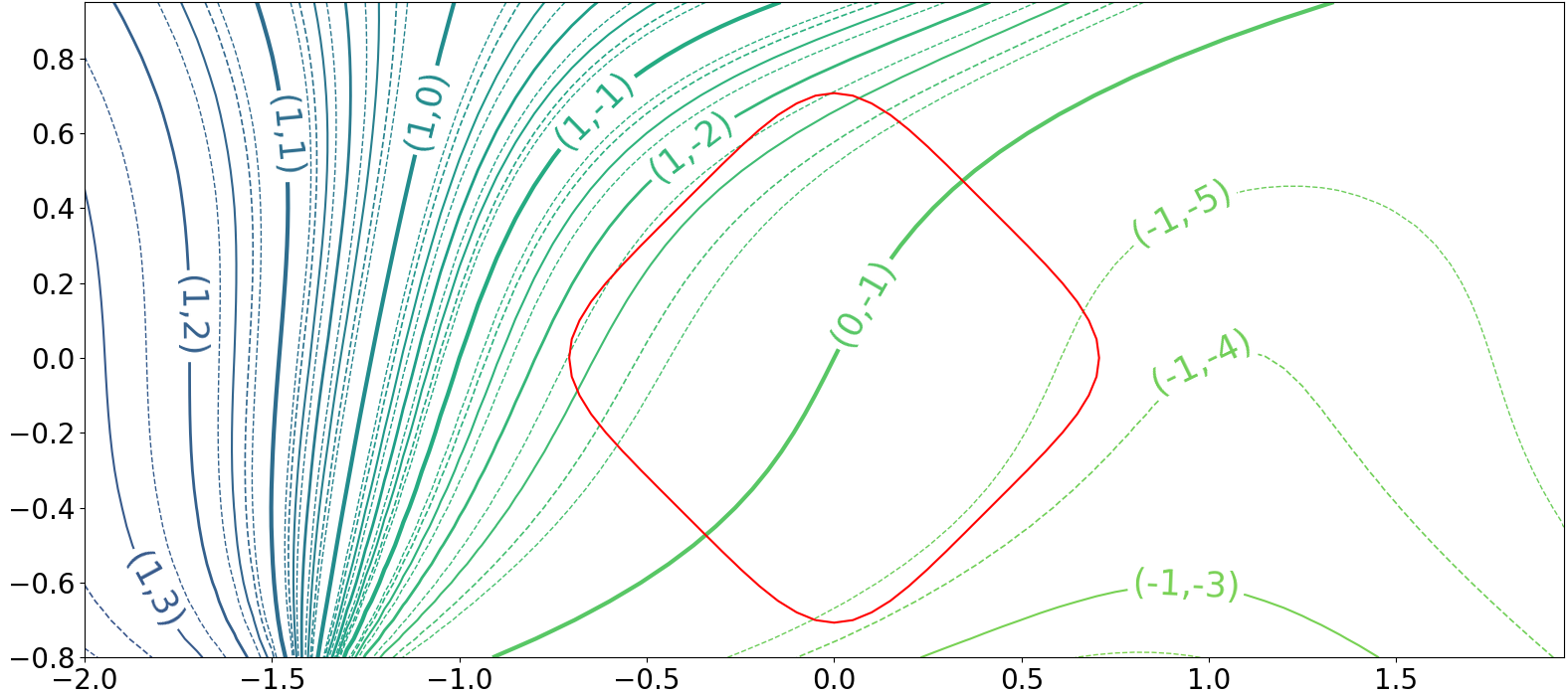}
	\centering
	\caption{The $(I_1, I_2)$-plane for Example~\ref{example.nonablnabl}. Inside the red curve (defined by $\det\Hess\: h(I)=0$), the Hessian $\Hess\: h(I)$ is indefinite. Several resonance lines $\Pi_j$ cross the red curve, but $\Gamma$ is empty and yields no intersection.} \label{fig.nonabnabl}
\end{figure}
\vspace{-2pt}

\begin{exa}
	$h(I) = I_1I_2$.
	\label{example.no.condC}
\end{exa}

This is a standard example of a non-convex Hamiltonian. The Hessian is constant and indefinite, with null directions $(1, 0)$ and $(0, 1)$. Moreover, $\Gamma = (\mathbb{R} \times \{ 0\} \cup \{0\} \times \mathbb{R}) \setminus \{(0, 0)\} $ is non-empty. However, $\Gamma = (\Pi_{0, 1} \cup \Pi_{1, 0}) \setminus \{(0, 0)\}$. So $\Gamma$ intersects some resonance lines non-transversally, and Condition \eqref{ConditionC} fails. In fact, this Hamiltonian is so degenerate that the steepness condition\footnote{A generic transversality condition, see e.g. \cite{N04}.} fails,
 and all stability estimates break down. Nevertheless, $h \in \mathcal{ID}(\R^2)$. Thus a generic deformation of $h$ admits meanders.

\begin{exa}
	$h(I) = I_1 + 3I_1 I_2 + 0.5 I_2 + I_2^2 + 0.6 I_2^3$.
	\label{example.yesnablnabl}
\end{exa}

\begin{figure}[ht!]
	\includegraphics[width=11cm]{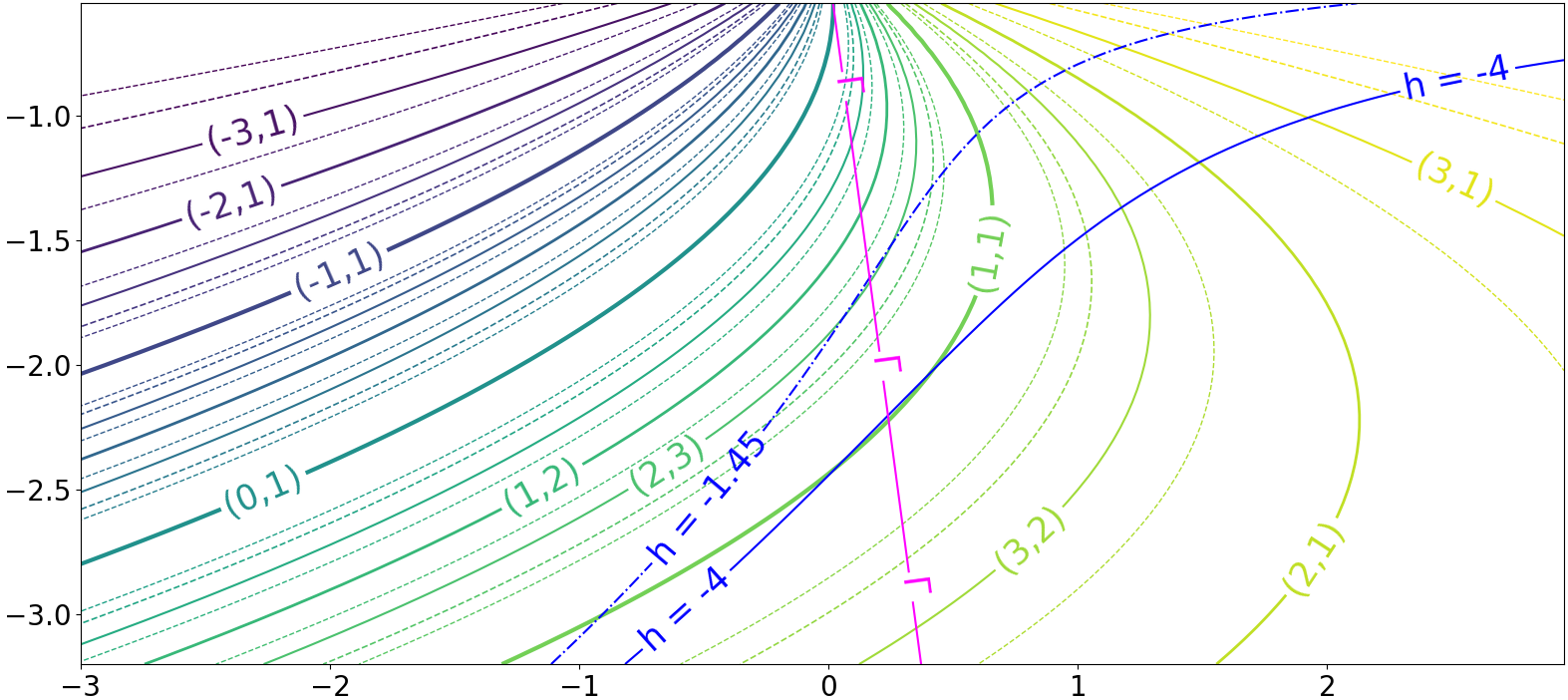}
	\centering
	\caption{The $(I_1, I_2)$-plane of Example \ref{example.yesnablnabl}. $\Pi_{1, 1}$ crosses $\Gamma$  transversely near the location $(0.24, -2.22)$ with energy close to $-4$;  $\Pi_{2,3}$ crosses $\Gamma$ transversely near $(0, -1.5)$ with energy close to $-1.45$.}
	\label{fig.yesnabnabl}
\end{figure}

This $h$ is generic isoenergetically degenerate and was chosen arbitrarily to illustrate the richness of the family. The submanifold $\Gamma$ crosses several resonance lines, as shown in Figure~\ref{fig.yesnabnabl}. Although $h$ is not quasi-convex, it is steep at the intersection points. Consequently, stability estimates remain valid, albeit with weaker constants.

\begin{figure}[ht!]
\begin{subfigure}{.5\textwidth}
  \centering
			         \includegraphics[width=5.5cm]{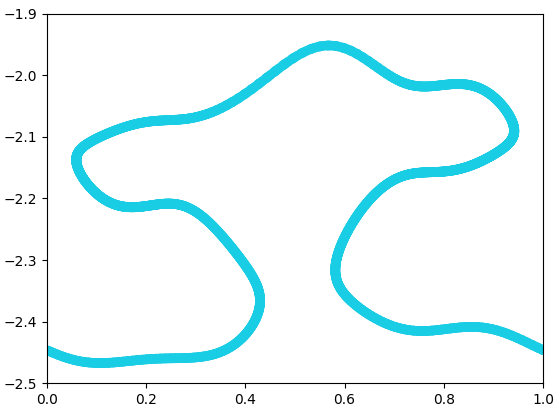}
                \medskip
\end{subfigure}%
\begin{subfigure}{.5\textwidth}
  \centering
  \includegraphics[width=5.5cm]{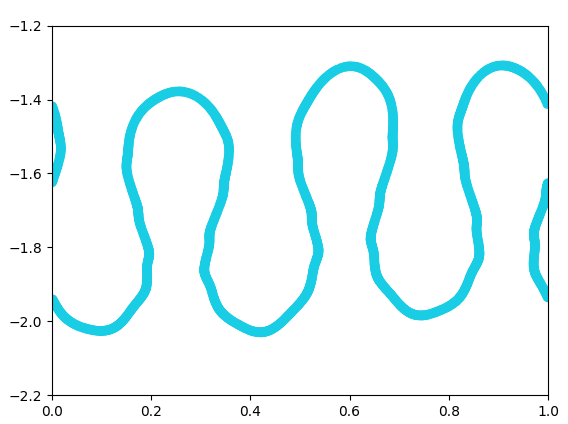}
                     \medskip
\end{subfigure}
\caption{Meandering invariant curves of Poincar\'e map for $h+p$, with energy $H = -4.08$ on the left, $H  =-1.45$ on the right.}
\label{ex.of.mush}
\end{figure}

Let $p = 0.2 \sin(2\pi \varphi_1) + 0.2 \sin(2\pi \varphi_2) + 0.02 \sin(4\pi \varphi_1 + 6\pi \varphi_2)$. Figure \ref{ex.of.mush} shows meanders of Poincar\'e map for $h+p$ corresponding to the intersections. 

\section{Structure of the proof}

A brief outline of the proof was given in Section ~\ref{sect.brief.pf}. We now explain each major stage in detail by presenting simplified versions of theorems. We also provide a roadmap of the proof, highlighting the order and role of each stage in Figure \ref{fig.planfirststep}.

\subsection{Constructing first step meanders}
The first step has two goals. The first is to construct a first-step meander. The second, relevant when a higher order meander is desired, is to carry out the first step of its construction. Both goals are achieved via the same scheme, as the higher-order meander is built upon the first-order one. So, the two cases diverge only at the final stage. We therefore start with describing the scheme of the first step meander construction.

\subsubsection*{Stage 1: Resonant Normal forms and approximation by one degree of freedom}

We start with a Hamiltonian with two degrees of freedom
\begin{equation}
	H(\varphi, I) = h(I) + p(\varphi, I), \qquad \varphi \in \T^2,\ \ \  I\in \R^2.
	\label{eq.hplusp}
\end{equation}
Hamiltonians with one degree of freedom are always integrable, with invariant curves given by energy levels. Therefore, a natural approach to studying invariant tori of a two-degree-of-freedom Hamiltonian is to reduce it to one degree of freedom. The resonant normal form provides precisely such a reduction, up to a  stretched exponential small error:

\begin{majthm}
	\label{maj.aver}
	If $\ \nabla h(I_*)\cdot (j_1,j_2)=0$, then there exists a symplectic coordinate change $L\circ \Psi$ such that in the new coordinates $H$ has a form
	\begin{equation}
		H(\varphi, I) = h(I) + g(\varphi_1, I) + f(\varphi, I), \qquad \varphi=(\varphi_1,\varphi_2),
	\end{equation}
	where $g(\varphi_1, I)$ is of order $\e$ and $f(\varphi, I)$ is  stretched exponential small in $\e$.
\end{majthm}

This coordinate change is a composition of $L$ and $\Psi$. $L$ is simply a linear coordinate change in $\SL(2,\Z)$ that reduces the case of arbitrary $(j_1,j_2)\in\Z^2\backslash\{0\}$ into $(1, 0)$. The second coordinate change $\Psi$ is more complicated, given by cohomological equation solving to get the resonant normal form, where \cite{P93} is applied.

\subsubsection*{Stage 2: Constructing a meandering channel with one degree of freedom}

Since $f$ is  stretched exponential small, we will first study the leading Hamiltonian $h(I)+g(\varphi_1, I)$. Since the system is independent of $\varphi_2$, the action $I_2$ is conserved and may be treated as a parameter. Consequently, $h(I) + g(\varphi_1, I)$ defines a one-degree-of-freedom Hamiltonian $\mathcal{H}_{I_2}(\varphi_1, I_1)$ parametrized by $I_2$.


\begin{figure}[ht!]
	\includegraphics[width=12cm]{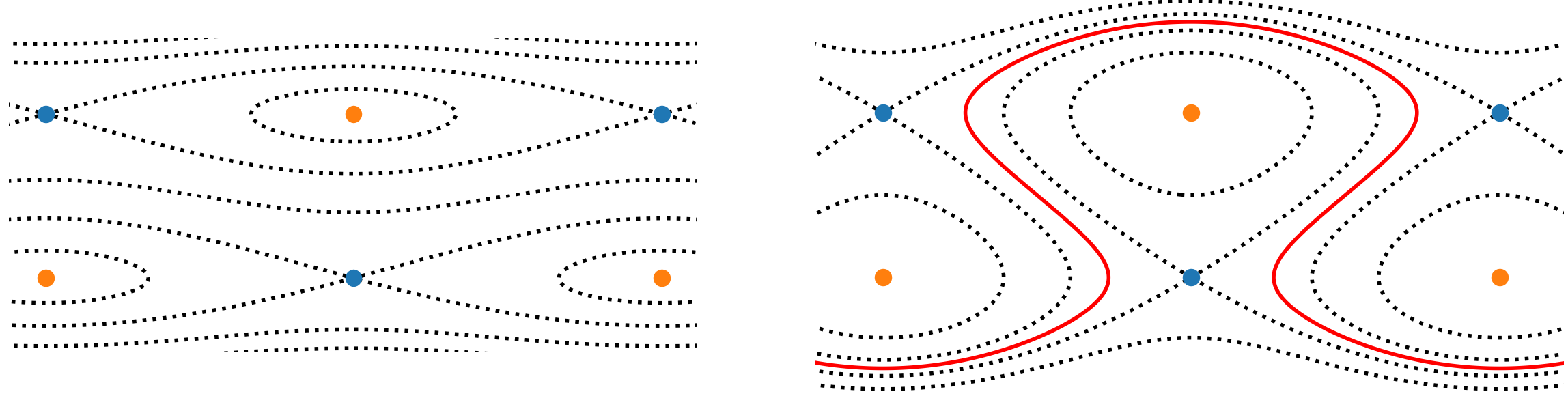}
	\centering
	\caption{Phase portrait of $h(I) + p(x, I) = I - I^3/3 + \varepsilon \cos 2\pi x$,  $\varepsilon = 0.2$ on the left and $\varepsilon = 1$ on the right. In the latter meandering invariant curves appear. }
	\label{fig.meander1d}
\end{figure}

We aim to show that the Hamiltonian $\cH$ has a meandering channel, where an example is shown in Figure \ref{fig.meander1d}. In fact, our non-degeneracy conditions guarantee the following result.

\begin{majthm}
Assume $h$ is generic isoenergetically degenerate and $g$ is $\|p\|_{r,s}^{1+\eta}$-non-degenerate. There exists a parameter $I_2^*$, such that the Hamiltonian
	\begin{equation}
		\cH_{I_2^*}(\varphi_1, I_1) = h(I_1, I_2^*) + g(\varphi_1, I_1, I_2^*)
	\end{equation}
	has a meandering channel.
	\label{maj.channel}
\end{majthm}



\subsubsection*{Stage 3: Action-angle estimates}

To extend the existence of invariant tori from $\mathcal{H}_{I_2^*}$ to $\mathcal{H}_{I_2^*} + f$, we will apply a KAM-type theorem. As KAM theorem requires the unperturbed Hamiltonian to depend on actions only, we pass to action-angle coordinates as follows.

\begin{majthm}
	\label{maj.actionangle}
	There exists a local action-angle coordinate change
$\Phi:(\theta, J) \rightarrow (\varphi, I)$,
so that locally we have \[(h+g)\circ \Phi(\theta,J)= h_{\two}(J).\]
	Moreover, its analytic properties (like the complex domain of definition and bounds on its image) are given in a quantitative way. 
\end{majthm}

We point out that analytic estimates are needed throughout. After all the coordinate changes above, the original system $h+p$ takes the form
\begin{equation}
	h_{\two}(J) + f \circ \Phi(\theta, J).
	\label{eq.hinactionangle}
\end{equation}
we require $f \circ \Phi$ to remain sufficiently small so that \eqref{eq.hinactionangle} can be regarded as a perturbation of $h_{\two}(J)$.


On one hand, the smallness of $f \circ \Phi$ allows us to apply a KAM theorem. In addition, KAM requires conditions on the integrable part $h_{\two}(J)$: namely, non-degeneracy and a Diophantine rotation number. 

On the other hand, if we seek for a higher order meander, \eqref{eq.hinactionangle} has the same structure as the original system \eqref{eq.hplusp}, with $h_{\two}(J)$ playing the role of $h$ and $f \circ \Phi$ that of $p$. As in the first step, we require $h_{\two}(J)$ to be generic isoenergetically degenerate, that is, to satisfy Conditions \eqref{ConditionA}--\eqref{ConditionE}. This is made precise by the following result.

\begin{majthm}
	The expressions in Conditions \eqref{ConditionA}--\eqref{ConditionE} for $h_{\two}(J)$ w.r.t. action-angle variables $(\theta, J)$ can be restated quantitatively in terms of expressions in the original variables $(\varphi, I)$. 
	\label{maj.actionangleconditio}
\end{majthm}


\subsubsection*{Stage 4: Choosing an invariant curve in the meandering channel}

We note that the conditions required for KAM and for the higher-order meander preparation are distinct. For instance, KAM requires a Diophantine rotation number, whereas Condition \eqref{ConditionE} requires it to be rational. Nevertheless, in both cases one can adjust the parameter $I_2$ to change the channel, and vary the choice of torus within the channel. Since choosing an invariant torus amounts to choosing $J$, we have the following result.

\begin{majthm}
Both of the following cases are realizable by appropriate choices of $I_2$ and $J$:
	\begin{enumerate}
		\item the rotation number is Diophantine and $h_{\two}$ is non-degenerate at $J$,
		\item $h_{\two}$ is generic isoenergetically degenerate at $J$.
	\end{enumerate}
	\label{maj.meanderselection}
\end{majthm}


The subsequent argument splits into two cases according to Theorem~\ref{maj.meanderselection}. Case (2) leads to the second step and the construction of a higher-order meander.  Case (1) allows us to apply a KAM theorem and conclude the existence of a first order meander. Applying the KAM theorem with the quantitative bounds of \cite{PoschelKAM}, we obtain the following.

\begin{majthm}
	If case (1) of Theorem \ref{maj.meanderselection} holds, then the Hamiltonian \eqref{eq.hinactionangle} admits an invariant torus close to that of $h_{\two}(J)$. 
	\label{maj.KAM}
\end{majthm}

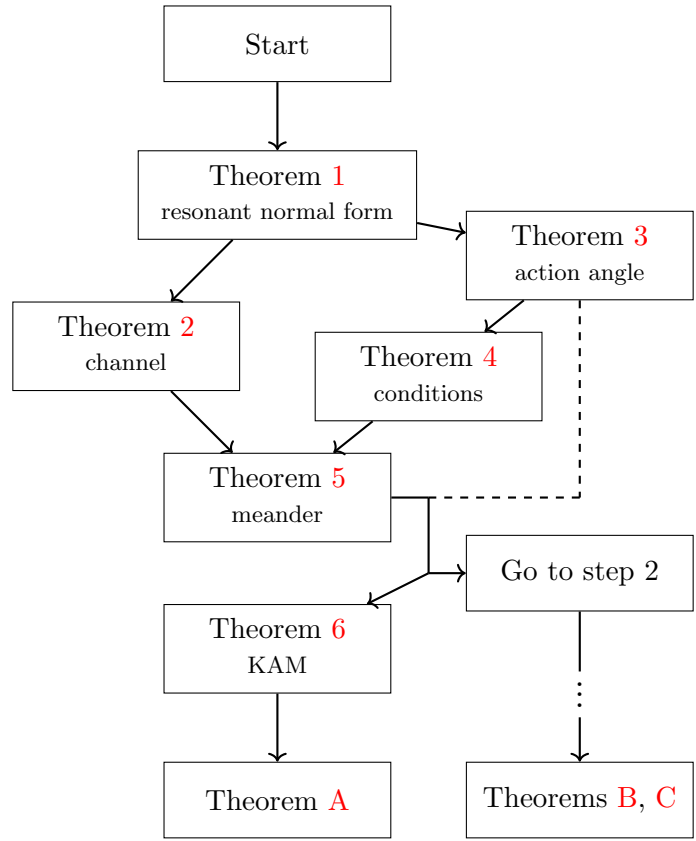
\begin{figure}[ht!]
	\begin{tikzpicture}[scale=2]
		
		\tikzstyle{box} = [rectangle, draw, minimum width=3cm, minimum height=1cm, align=center]
		
		\node[box] (top) at (0, 0) { \begin{tabular}{c} Theorem \ref{maj.aver} \\ \Small resonant normal form \end{tabular}};
		\node[box] (middle) at (-1, -1) { \begin{tabular}{c} Theorem \ref{maj.channel} \\ \Small channel \end{tabular}};
		\node[box] (bottom) at (0, -2) { \begin{tabular}{c} Theorem \ref{maj.meanderselection} \\ \Small meander \end{tabular}};
		\node[box] (action) at (2, -0.4) { \begin{tabular}{c} Theorem \ref{maj.actionangle} \\ \Small action angle \end{tabular}};
		\node[box] (action2) at (1, -1.2) {\begin{tabular}{c} Theorem \ref{maj.actionangleconditio} \\ \Small conditions \end{tabular}};
		\coordinate (help) at (2, -2);
		\coordinate (help2) at (1, -2);
		\coordinate (help3) at (1, -2.5);
		\node[box] (action3) at (0, -3) {\begin{tabular}{c} Theorem \ref{maj.KAM} \\ \Small KAM \end{tabular}};
		\node[box] (action4) at (2, -2.5) {Go to step $2$};
		\node[box] (action5) at (0, -4) {Theorem \ref{thm.A}};
		\node[box] (action6) at (2, -4) {Theorems \ref{thm.B}, \ref{thm.C}};
		\node[box] (action7) at (0, 1) {Start};
		
		\node (action8) at (2, -3.25){\Large $\mathbf \vdots$};
		\coordinate (help4) at (2, -3.15);
		\draw[->, thick] (top) -- (middle);
		\draw[->, thick] (middle) -- (bottom);
		\draw[->, thick] (top) -- (action);
		\draw[->, thick] (action) -- (action2);
		\draw[->, thick] (action2) -- (bottom);
		\draw [dashed, thick] (action) -- (help);
		\draw [dashed, thick] (help) -- (help2);
		\draw [thick] (bottom) -- (help2);
		\draw [thick] (help2) -- (help3);
		\draw[->, thick] (help3) -- (action3);
		\draw[->, thick] (help3) -- (action4);
		\draw[->, thick] (action3) -- (action5);
		\draw[->, thick] (action7) -- (top);
		\draw[thick] (action4) -- (help4);
		\draw[->, thick] (action8) -- (action6);
	\end{tikzpicture}
	\caption{A plan of the proof on the first step. The solid arrows indicate a major dependency, the dashed -- minor. At the end we can either finish with the first meander and apply KAM with Theorem \ref{maj.KAM} or go to the second step.}
	\label{fig.planfirststep}
\end{figure}

Finally, we transfer the invariant torus given by Theorem~\ref{maj.KAM} back to the original coordinates. In this way we deduce a meander and complete the proof of Theorem \ref{thm.A}.The plan of the first step is summarized in Figure \ref{fig.planfirststep}.

\subsection{Constructing second step meander}
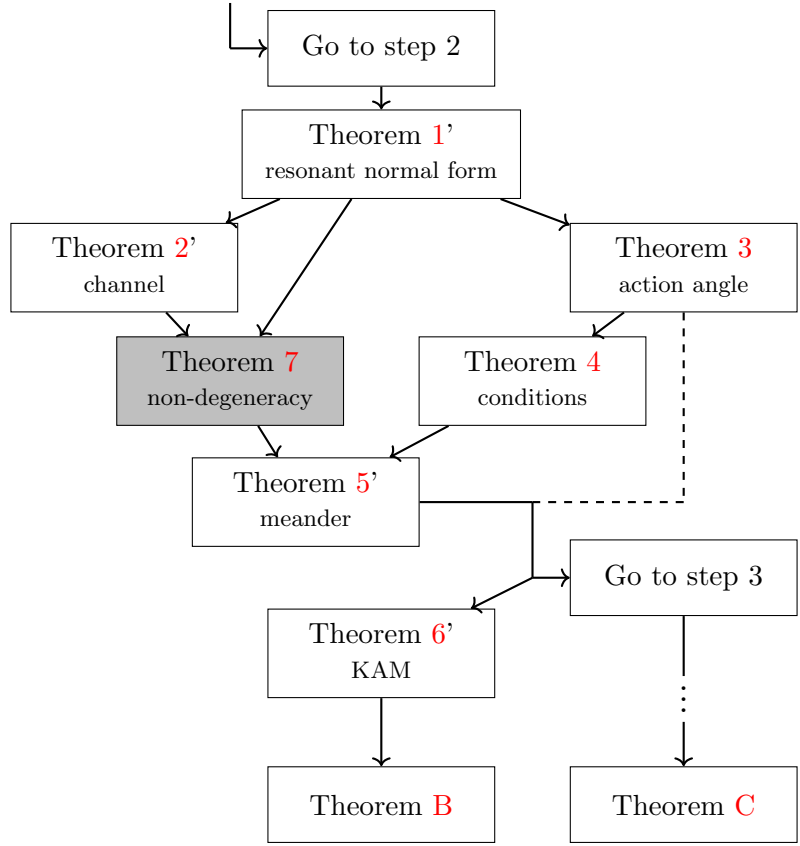
\begin{figure}[ht!]
	\begin{tikzpicture}[scale=2]
		
		\tikzstyle{box} = [rectangle, draw, minimum width=3cm, minimum height=1cm, align=center]
		
		\node[box] (top) at (0, 0.3) { \begin{tabular}{c} Theorem \ref{maj.aver}' \\ \Small resonant normal form \end{tabular}};
		\node[box] (left) at (-1.7, -0.45) { \begin{tabular}{c} Theorem \ref{maj.channel}' \\ \Small channel \end{tabular}};
		\node[box, fill=lightgray] (middle) at (-1, -1.2) { \begin{tabular}{c} Theorem \ref{maj.nondegenerate} \\ \Small non-degeneracy \end{tabular}};
		\node[box] (bottom) at (-0.5, -2) { \begin{tabular}{c}  Theorem \ref{maj.meanderselection}' \\ \Small meander \end{tabular}};
		\node[box] (action) at (2, -0.45) { \begin{tabular}{c}  Theorem \ref{maj.actionangle} \\ \Small action angle \end{tabular}};
		\node[box] (action2) at (1, -1.2) {\begin{tabular}{c} Theorem \ref{maj.actionangleconditio} \\ \Small conditions \end{tabular}};
		\coordinate (help) at (2, -2);
		\coordinate (help2) at (1, -2);
		\coordinate (help3) at (1, -2.5);
		\node[box] (action3) at (0, -3) {\begin{tabular}{c} Theorem \ref{maj.KAM}' \\ \Small KAM \end{tabular}};
		\node[box] (action4) at (2, -2.5) {Go to step $3$};
		\node[box] (action5) at (0, -4) {Theorem \ref{thm.B}};
		\node[box] (action6) at (2, -4) {Theorem \ref{thm.C}};
		\node[box] (action7) at (0, 1) {Go to step $2$};
		\coordinate (help7) at (-1, 1);
		\coordinate (help8) at (-1, 1.3);

		\node (action8) at (2, -3.25){\Large $\mathbf \vdots$};
		\coordinate (help4) at (2, -3.15);
		\draw[->, thick] (top) -- (middle);
		\draw[->, thick] (top) -- (left);
		\draw[->, thick] (left) -- (middle);
		\draw[->, thick] (middle) -- (bottom);
		\draw[->, thick] (top) -- (action);
		\draw[->, thick] (action) -- (action2);
		\draw[->, thick] (action2) -- (bottom);
		\draw [dashed, thick] (action) -- (help);
		\draw [dashed, thick] (help) -- (help2);
		\draw [thick] (bottom) -- (help2);
		\draw [thick] (help2) -- (help3);
		\draw[->, thick] (help3) -- (action3);
		\draw[->, thick] (help3) -- (action4);
		\draw[thick] (help8) -- (help7);
		\draw[->, thick] (help7) -- (action7);
		\draw[->, thick] (action3) -- (action5);
		\draw[->, thick] (action7) -- (top);
		\draw[thick] (action4) -- (help4);
		\draw[->, thick] (action8) -- (action6);
	\end{tikzpicture}
	\caption{A plan of the proof on the second step. Gray boxes represent theorems that don't appear on the first step. These theorems are the most difficult part of the second step.}
	\label{fig.secondstepplan}
\end{figure}

We now outline the proof of the second step. Later steps follow the same plan and differ only in estimates. The goal is Theorem~\ref{thm.B}, and the plan is summarized in Figure~\ref{fig.secondstepplan}, where analogues of the first-step theorems are marked with a prime.

By case (2) of Theorem~\ref{maj.meanderselection}, the second step begins with $h_{\two}(J) + f \circ \Phi(\theta, J)$, where $h_{\two}$ is generically isoenergetically degenerate and $f\circ \Phi(\theta, J)$ is small. The overall framework is thus analogous to the first step.

The key difficulty, however, is that Theorem~\ref{maj.meanderselection} provides no information on $f\circ \Phi$ beyond smallness. In particular, it does not supply the non-degeneracy conditions \eqref{ConditionF}--\eqref{ConditionH}, which are essential for meanders to appear. This is the main technical obstacle of the second step.

\subsubsection{Introducing corrector $\mu$}

To address this difficulty, we allow a small modification $p$. Replacing $p$ by $p+\mu$ introduces a correction to $f$
 $$f \rightarrow f + \mu \circ L \circ \Psi,$$
 giving us some control over $f$. Since $f$ arises from the resonant normal form (Theorem~\ref{maj.aver}), to guarantee that the corrector $\mu$ does not average out, $\mu $ must be chosen small enough not to exceed the upper bound on $f$ given there.

However, the irregularity of $\Psi$ forces $\mu$ to be smaller than any polynomial of $\|f\|$. This correction is therefore too small to directly establish Conditions~\eqref{ConditionF}--\eqref{ConditionH}. Instead, we work with more delicate non-degeneracy conditions of types I--IV.

As in the first step, a parallel resonant normal form Theorem~\ref{maj.aver}  yields  $g_{\two}$ and $f_{\two}$, and an analogue of Theorem~\ref{maj.channel} provides a channel for some parameter $I_2$. However, without any non-degeneracy assumption on $g_{\two}$, this channel may be degenerate. It is precisely where the corrector $\mu$ plays its role.

\subsubsection{Resolving degeneracies with $\mu$}
We establish the following:
\begin{majthm}
There exists $\mu$ such that, for some value of $I_2$, the Hamiltonian $h_{\two} + g_{\two}$ admits a non-degenerate meandering channel.
	\label{maj.nondegenerate}
\end{majthm}

Constructing an analytic $\mu$ explicitly turns out to be difficult. Instead, we work with a smooth approximation and recover real analyticity afterwards. We consider a parametrized family of functions and estimate the measure of the set of parameters for which the resulting channel is degenerate. By a careful choice of family, this degenerate set has sufficiently small measure, yielding the genericity condition on $p$ in the second step.

In particular, we point out that by pulling back the smooth functions to the original coordinates and approximating by analytic functions on $W_s\T^2\times V_rD$, we manage to deform $p$ without loss of analyticity. 

The remainder of the second step parallels the first. Theorems~\ref{maj.actionangle} and~\ref{maj.actionangleconditio} apply directly. Once the channel is non-degenerate, an analogue of Theorem~\ref{maj.meanderselection} guarantees that both the KAM case and the generically isoenergetically degenerate case are realizable for certain invariant torus within the channel. In the former case, we conclude the second step via a KAM theorem parallel to Theorem~\ref{maj.KAM}; in the latter, we proceed to the third step.

\subsection{Constructing general step meander}
To construct a Hamiltonian with meanders of arbitrary or infinite order, we iterate the procedure from the second step. All subsequent steps follow the same plan as the second step and Theorem~\ref{thm.C} is established by induction. 

A key feature of the iteration is that the deformation $\mu_{(n)}$
at step $n$ is extremely small. To be precise, its size is bounded by an exponential tower in $1/\varepsilon$ with $n-1$ exponentials. To construct a meander of infinite order, the smallness condition on $\varepsilon$ must not accumulate additional constraints at each step. In fact, we prove that all inequalities involving exponential towers hold uniformly in $n$. And the construction of deformation by analytic approximation as in second step preserves the region of analyticity. 
This completes the outline of the proof.

\subsection*{Acknowledgement}
This work is partially supported by ERC Advanced Grant SPERIG ($\# 885707$). The authors are grateful to Marie-Claude Arnaud, Abed Bounemoura, Luigi Chierchia, Kostiantyn Drach, Holger Dullin, Anna Florio, Rapha\"el Krikorian, Jean-Pierre Marco, Frank Trujillo, Ke Zhang for discussion and reference. Y.P. would like to thank the hospitality of Institute Mittag-Leffler and University of Maryland.

\printbibliography

\end{document}